\documentclass[12pt]{amsart}
\usepackage[dvips]{graphicx}

\newtheorem{theorem}{Theorem}

\theoremstyle{definition}

\theoremstyle{remark}
\newtheorem{remark}[theorem]{Remark}
\numberwithin{equation}{section}

\usepackage{amsmath,amssymb,txfonts}
\usepackage{fancybox} 
\usepackage{ulem}
\usepackage{wrapfig}
\usepackage{comment}
\usepackage{amsthm,color,bm}

\def\re{\mathbf{R}}
\def\N{\mathbf{N}}
\def\Z{\mathbf{Z}}
\def\tc{\textcolor}
\def\({\left(}
\def\){\right)}
\def\pd{\partial}

\def\ep{\varepsilon}
\def\w{\omega}
\def\la{\lambda}

\def\la{\lambda}

\def\la{\lambda}

\begin{document}
\title[Two limits on Hardy and Sobolev inequalities]{Two limits on Hardy and Sobolev inequalities}

\author{Megumi Sano}
\address{Laboratory of Mathematics, Graduate School of Engineering,
Hiroshima University\\
Higashi-Hiroshima, 739-8527, Japan}
%\curraddr{}
\email{smegumi@hiroshima-u.ac.jp}
%\thanks{}

%\subjclass[2010]{Primary 35A23; Secondary 26D10.}

%\keywords{Hardy inequality, remainder term, virtual extremal.}
\date{\today}

\dedicatory{}

\begin{abstract}
It is known that classical Hardy and Sobolev inequalities hold when the exponent  $p$ and the dimension $N$ satisfy $p < N < \infty$. 
In this note, we consider two {\it limits} of Hardy and Sobolev inequalities as $p \nearrow N$ and $N \nearrow \infty$ in some sense.
\end{abstract}

\maketitle

%%
%% Start line numbering here if you want
%%
% \linenumbers

\tableofcontents

%%%%%%%%%%%%%%%%%%%%%%%%%%%%%%%%%%%%%%%%%%%%%%%%%%%%%%%%%%%%%%%%%%%%%%%%%%%%%%%%%%%%%%%%%%%%%%%%%%%%%%%%%%%%%%%
%% この下から本文を編集してください．
%%%%%%%%%%%%%%%%%%%%%%%%%%%%%%%%%%%%%%%%%%%%%%%%%%%%%%%%%%%%%%%%%%%%%%%%%%%%%%%%%%%%%%%%%%%%%%%%%%%%%%%%%%%%%%%
\section{Introduction : Hardy and Sobolev inequalities}

%%%%%%%%%%%%%%%%%%%%%%%%%%%%%%%%%%%%%%%%%%%%%%%%%%%%%%%%%%%%%%%%
%
%\section{導入 : Hardy不等式とSobolev不等式}
%
%%%%%%%%%%%%%%%%%%%%%%%%%%%%%%%%%%%%%%%%%%%%%%%%%%%%%%%%%%%%%%%%

Let $\Omega \subset \re^N$ be a domain and $0 \in \Omega$. 
If the exponent $p \ge 1$ and the dimension $N \ge 2$ satisfy $p < N$, then the Hardy inequality (\ref{H}) and the Sobolev inequality (\ref{S}) hold for any $u \in W_0^{1,p} (\Omega )$, where $W_0^{1,p} (\Omega )$ is a completion of $C_c^{\infty}(\Omega )$ with respect to $\| \nabla \cdot \|_{L^p(\Omega)}$.
\begin{align}\label{H}
\( \dfrac{N-p}{p} \)^p \int_{\Omega} \dfrac{|u(x)|^p}{|x|^p} dx &\le \int_{\Omega} | \nabla u(x) |^p \, dx,
\end{align}
\begin{align} 
\label{S}
S_{N,p} \( \int_{\Omega} |u(x)|^{p^*} \,dx \)^{\frac{p}{p^*}} &\le \int_{\Omega} | \nabla u(x) |^p \, dx, \\
\text{where}\,\,p^* = \frac{Np}{N-p}, \,&S_{N,p} = \pi^{\frac{p}{2}} N \( \frac{N-p}{p-1} \)^{p-1} \( \frac{\Gamma (\frac{N}{p}) \Gamma (N+ 1 -\frac{N}{p})}{\Gamma (N)  \Gamma(1+\frac{N}{2})} \)^{\frac{p}{N}}. \notag
\end{align}
These two inequalities appear in analyzing existence, non-existence, and stability of solution to nonlinear partial differential equations and so on. And their best constants and their attainability are well-studied (ref. \cite{Au}, \cite{T}, \cite{L}, \cite{BV}, \cite{BG}, \cite{GR}, \cite{RS} etc.). 
Not only that, the Sobolev inequality (\ref{S}) denotes an embedding of the subcritical Sobolev space : $W_0^{1, p} \hookrightarrow L^{p^*}$, and 
the Hardy inequality (\ref{H}) denotes an embedding : $W_0^{1, p} \hookrightarrow L^{p^*, p} (\subsetneq L^{p^*})$. 
Therefore these two inequalities are fundamental and important. 

The Hardy inequality (\ref{H}) and the Sobolev inequality (\ref{S}) hold only when the exponent $p$ and the dimension $N$ satisfy
$$p < N \,\,(\text{the subcritical case}).$$
Indeed, we observe that if taking limits as $p \nearrow N$, two best constants $(\frac{N-p}{p})^p, S_{N,p}$ go to zero and two integrals $\| u \|_{L^p( |x|^{-p}dx)}, \| u\|_{L^{p^*}}$ for a suitable function $u$ diverge, here $\frac{1}{0} = \infty$. 
Therefore we see that two inequalities break down due to the indeterminate forms: $0 \times \infty$. 
However we know well about a sequence of real numbers, it is possible to exist its limit even if it is an indeterminate form. 
In the same spirit, by making two quantities compete with each other, of which one goes to zero, and the other diverges by taking a limit of an exponent in functional inequality, 
{\bf can we get some ``limit'' of functional inequality?} 
And {\bf can we get some ``infinite dimensional form'' of functional inequality?} 

In this note, we explain limiting procedures for the Hardy inequality (\ref{H}) and the Sobolev inequality (\ref{S}) as $p \nearrow N$ and $N \nearrow \infty$.\\

\underline{\bf Notation}
\begin{itemize}
\item $B_R^N = \{ x \in \re^N \,\, : \,\, |x| < R \,\}$.
\item $\w_{N-1} = \frac{N \pi^{\frac{N}{2}}}{\Gamma \( 1+ \frac{N}{2} \)}$：the area of the unit sphere $\mathbb{S}^{N-1} \subset \re^N$.
\item $\Gamma (t)= \int_0^\infty x^{t-1} e^{-x} \,dx$：the Gamma function.
\item $u^* (t)= \inf \{ \tau >0\,\,:\,\, |\{ x \in \re^N : |u(x)| > \tau \}| \le t \}$：the rearrangement of $u$.
\item $L^{p,q}(\log L)^r = \left\{ \, u: \Omega \to \re\,\,{\rm measurable} \,\,: \,\, \| u\|_{L^{p,q} (\log L)^r} < \infty  \, \right\}$：the Lorentz-Zygmund space.
\end{itemize}
\begin{align*}
\| u \|_{L^{p,q}(\log L)^r} &= 
\begin{cases}
\( \int_0^{|\Omega|} s^{\frac{q}{p}-1} \( e+ |\log s| \)^{rq} u^* (s)^q \,ds \)^{\frac{1}{q}} \quad &\text{if} \,\, 1\le q <\infty, \\
\sup_{0< s < |\Omega|} s^{\frac{1}{p}} \( e+ |\log s| \)^r u^* (s)   &\text{if} \,\, q = \infty.
\end{cases}
\end{align*}
Note that $\| \cdot \|_{L^{p,q}(\log L)^r}$ is not norm. That is a quasi norm. 
Moreover $L^{p,q}(\log L)^0$ is the Lorentz space $L^{p,q}$, $L^{\infty, \infty}(\log L)^r$ is the Zygmund space $Z^{-r}$, and Zygmund space $Z^{-r}$ coincides with the Orlicz space $L_{e^{\,|u|^{-1/r}}} (= {\rm ExpL}^{-\frac{1}{r}})$ with the  Young function $\Phi (t)= e^{|t|^{-1/r}}-1$ (see \cite{BR} p.15 Theorem D (c)).

%
%
%%%%%%%%%%%%%%%%%%%%%%%%%%%%%%%%%%%%%%%%%%%%%%%%%%%%%%%%%%%%%%%%
%
%\section{Indirect limiting procedure}
%
%%%%%%%%%%%%%%%%%%%%%%%%%%%%%%%%%%%%%%%%%%%%%%%%%%%%%%%%%%%%%%%%
%
\section{Indirect limiting procedure}

\subsection{$p \nearrow N$}
\tc{white}{a}\\
Fix the dimension $N\ge 2$. 
We consider some {\it limit} of the Hardy inequality (\ref{H}) and the Sobolev inequality (\ref{S}) as $p \nearrow N$. 
There is no general theory of getting a limiting form of the Hardy inequality (\ref{H}) and the Sobolev inequality (\ref{S}). 
Therefore we have to think out {\bf how to take a limit} corresponding to each case. 
However, for these two inequalities, roughly speaking, we {\it divide} a subcritical information, and we {\it sum (or combine)} them, then we can get a critical (limit) information. 

\subsubsection{The Sobolev inequality}

\begin{theorem}\label{IL S} (\cite{Yu, P, Tru})
Let $|\Omega | < \infty$. 
We obtain the following non-sharp Trudinger-Moser inequality (\ref{TM}) as a limit of the Sobolev inequality (\ref{S}) as $p \nearrow N$.
\begin{align}\label{TM}
\int_{\Omega} {\rm exp} \left[ \alpha \( \frac{|u(x)|}{\| \nabla u \|_N} \)^{\frac{N}{N-1}} \right] \, dx \le C \,| \Omega | \quad \text{for small} \,\, \alpha >0. 
\end{align}
\end{theorem}

\begin{remark}
Here, non-sharp means it does not have its optimal exponent and its best constant. 
It is known that the optimal exponent $\alpha$ of the  Trudinger-Moser inequality (\ref{TM}) is $N \w_{N-1}^{\frac{1}{N-1}}$ (ref. \cite{M}). 
However in Theorem \ref{IL S}, we can not obtain information of the optimality. 
Unfortunately, we do not know the exact value of the best constant $C$ even now.
\end{remark}

Theorem \ref{IL S} was shown by Yudovich \cite{Yu}, Pohozaev \cite{P}, Trudinger \cite{Tru}. 
Here we prove it by the proof of \cite{CR} Theorem 1.7., which uses information  of the decay speed of Sobolev's best constant $S_{N,p}$ as $p \nearrow N$.  \\

%%%%%%%%%%%%%%%%%%%%%%%%%%%%%%%%%%%%%%

\noindent
{\it Proof of Theorem \ref{IL S}}：\,\,
The decay order of Sobolev's best constant $S_{N,p}$ as $p \nearrow N$ is as follows. 
$$
S_{N,p} \sim (N-p)^{p-1} \quad (p \nearrow N).
$$
For $q \in (N, \infty)$, set $p=p(q) = \frac{Nq}{N+q}$. 
Thenwe see that $p \nearrow N \iff q \nearrow \infty$ and $p^* = q$. 
By the Sobolev inequality (\ref{S}), we have
\begin{align*}
\| u\|_q = \| u\|_{p^*} \le S_{N,p}^{-\frac{1}{p}} \| \nabla u\|_p 
&\le C (N-p)^{-\frac{p-1}{p}} |\Omega |^{\frac{1}{p} - \frac{1}{N}} \| \nabla u\|_N \\
&=C \( \frac{N+q}{N^2} \)^{\frac{N-1}{N} - \frac{1}{q}} |\Omega|^{\frac{1}{q}} \, \| \nabla u\|_N.
\end{align*}
Thus for any $q \in (N, \infty)$, we have
\begin{align}\label{q order}
\| u\|^q_q \le C^q \,q^{\frac{N-1}{N}q - \frac{1}{q}} |\Omega| \,\| \nabla u \|_N^q.
\end{align}
Applying the inequality (\ref{q order}) for $q=\frac{N}{N-1} k$ implies that
\begin{align*}
\int_{\Omega} {\rm exp} \left[ \alpha \( \frac{|u(x)|}{\| \nabla u \|_N} \)^{\frac{N}{N-1}} \right] \, dx
&= \sum_{k=0}^\infty \int_{\Omega} \frac{\alpha^k}{k!} \( \frac{|u(x)|}{\| \nabla u \|_N} \)^{\frac{N}{N-1}k} \,dx\\
&\le C + |\Omega|\sum_{k=M}^\infty \frac{\alpha^k C^{\frac{N}{N-1} k} \, \(  \frac{Nk}{N-1} \)^{k-1} }{k!},
\end{align*}
where $M \gg 1$. 
By Stirling's formula $k! \sim \sqrt{2\pi k} \,k^k e^{-k}\,(k \to \infty)$, we see that the right-hand side of the above inequality does not diverge if
$$
\left| \frac{(\alpha C e )^k}{\sqrt{2\pi k}} \right|^{\frac{1}{k}} < 1.
$$
Therefore the inequality (\ref{TM}) holds for small $\alpha$. 
\qed

%%%%%%%%%%%%%%%%%%%%%%%%%%%%%%%%%%%%%%%%

\subsubsection{The Hardy inequality}
\tc{white}{a}\\
In a limiting case $p=N$ of the Hardy inequality, the following inequality which is called the critical Hardy inequality is known.
\begin{align}\label{CH}
\( \frac{N-1}{N} \)^N \int_{B^N_R} \frac{| u |^N}{|x|^N \( \log \frac{aR}{|x|} \)^N} \,dx \le \int_{B_R^N} | \nabla u |^N \, dx \quad (a \ge1).
\end{align}
The critical Hardy inequality (\ref{CH}) was founded by Leray \cite{Le}. However, unlike the Sobolev inequality, the critical Hardy inequality was not drived as a {\it limit} of the Hardy inequality (\ref{H}) as $p \nearrow N$. Therefore, it is unclear, at least for me, that {\bf why the inequality (\ref{CH}) is called the critical Hardy inequality, and why the logarithmic function appears in the Hardy potential in a limiting form}.
To resolve it, we derive the logarithmic function via some limiting procedure for the Hardy inequality (\ref{H}) as $p \nearrow N$.

\begin{theorem}\label{IL H}(\cite{SS})
We obtain the following non-sharp critical Hardy inequality (\ref{NS CH}) as a limit of the Hardy inequality (\ref{H}) as $p \nearrow N$.
\begin{align}\label{NS CH}
C_{\beta, a} \int_{B^N_R} \frac{| u |^N}{|x|^N \( \log \frac{aR}{|x|} \)^\beta} \,dx \le \int_{B_R^N} | \nabla u |^N \, dx \quad (a > 1, \beta \gg 1).
\end{align}
\end{theorem}

\begin{remark}
$\beta \gg 1$ in Theorem \ref{IL H} is corresponding to $\alpha \ll 1$ in Theorem \ref{IL S}. Furthermore, in Theorem \ref{IL H}, since we forcus on  Hardy's best constant $(\frac{N-p}{p})^p$ and the singurality of the Hardy potential $|x|^{-p}$ at the origin only, we do not obtain the non-sharp inequality (\ref{NS CH}) with $a=1$ which has the boundary singurality. 
However, if we consider taking a {\it limit} of the Poincar\'e inequality in a domain $\Omega$ as $| \Omega | \searrow 0$ by using an information that  Poincar\'e's best constant $\la (\Omega)$ goes to infinity, it is possible to obtain the non-sharp inequality (\ref{NS CH}) with $a=1$. We omit here.
\end{remark}

%%%%%%%%%%%%%%%%%%%%%%%%%%%%%%%%%%%%%%
Without loss of generality, we assume $R=1$. 
Before the proof of Theorem \ref{IL H}, we prepare taking a limit of the Hardy inequality (\ref{H}) as $p \nearrow N$. 

Set $p_k=N-\frac{1}{k}$ for $k \in \N$. Then $k \nearrow \infty \iff p_k \nearrow N$. Since Hardy's best constant for the exponent $p_k\, (< N)$ is $(\frac{N-p_k}{p_k})^{p_k} \sim k^{-N}$, it goes to zero as $k \nearrow \infty$. 
On the other hand, since the Hardy potential $|x|^{-p_k}$ goes to $|x|^{-N} \not\in L^1(B_\ep)$ as $k \nearrow \infty$, the integral $\int_{B_1} \frac{|u|^{p_k}}{|x|^{p_k}} \,dx$ goes to infinity as $k \nearrow \infty$. 
In order to measure the speed of its divergence and make the integral compete Hardy's best constant, we consider the followings.

Let $f \in C^1 (0, \infty)$ be a monotone-decreasing function with $\lim_{t \to +\infty} f(t) =0$, and $\{ \phi_k \}_{k \in \Z} \subset  C_c^{\infty}(\re^N \setminus \{ 0\})$ be a sequence of radial functions which satisfy the followings.
\begin{align*}
&(i) \,\sum_{k=-\infty}^{+\infty} \phi_k (x)^N =1, \,0 \le \phi_k (x) \le 1 \,\, \( \forall x \in \re^N \setminus \{ 0\} \), \\
&(ii) \,\,\text{supp}\,\phi_k \subset B_{f(k)} \setminus B_{f(k+2)}.
\end{align*}
For a radial function $u \in C_c^1 (B_1)$, we set $u_k = u \,\phi_k, A_k= {\rm supp}\,u_k \subset B_1 \cap \( B_{f(k)} \setminus B_{f(k+2)} \)$. 

We can divide $\re^N \setminus \{ 0\}$ by balls $B_{f(k)}$. 
Whether we can obtain a limit of the Hardy inequality (\ref{H}) depends on $f$ which decides how to devide the domain $B_1$. 
In order to obtain a limit of the Hardy inequality (\ref{H}) by this limiting procedure, the left-hand side of the Hardy inequality (\ref{H}) with the exponent $p_k$ and the function $u_k$ must not be vanishing as $k \to \infty$. We shall determine such $f$.

Since $k \le f^{-1}(|x|) \le k+2$ for $x \in A_k$, the left-hand side of the Hardy inequality (\ref{H}) with $p_k$ and $u_k$ can be estimated as follows.
\begin{align}\label{H_p_k}
&\( \frac{N-p_k}{p_k} \)^{p_k} \int_{A_k} \frac{|u_k|^{p_k}}{|x|^{p_k}} dx 
= p_k^{-p_k} \int_{A_k} \( \, \frac{|u_k (x)|}{|x| \,k} \, \)^{N-\frac{1}{k}} dx \notag \\
&\ge C \int_{A_k} \frac{|u_k (x)|^N}{|x|^N \( f^{-1}(|x|) \)^N} \( \, \frac{|x|\, k}{|u_k (x)|} \, \)^{\frac{1}{k}} dx \notag \\
&\ge C \, \| \nabla u_k \|_{L^N(A_k)}^{-\frac{1}{k}} \int_{A_k} \frac{|u_k (x)|^N}{|x|^N \( f^{-1}(|x|) \)^N} \(  f(k+2) \( \log \frac{f(k)}{f(k+2)} \)^{-\frac{N-1}{N}} \, \)^{\frac{1}{k}} dx,
\end{align}
where the second inequality comes from the pointwise estimate (Radial lemma) for the radial function $u_k$\,: $|u_k (x)| \le \| \nabla u_k \|_N \( \log \frac{f(k)}{|x|} \)^{\frac{N-1}{N}} (x \in A_k)$.  
Therefore, if for any $k \in \N$ the function $f$ satisfies
\begin{align}\label{f}
\(  f(k+2) \( \log \frac{f(k)}{f(k+2)} \)^{-\frac{N-1}{N}} \, \)^{\frac{1}{k}} \ge C >0,
\end{align}
then the information on the left-hand side of the classical Hardy inequality (\ref{H}) is not vanishing in this limiting procedure. From (\ref{f}) and l'H\^opital's rule, we have an ordinary differential inequality for $f$ as follows:
\begin{align*}
\frac{d}{dt} f (t) \ge -C f (t)
\end{align*}
whose solution satisfies $f(t) \ge e^{-Ct}$. Thus $f^{-1}(t) \ge \frac{1}{C} \log \frac{1}{t}$. We belive that the above caluculation and consideration give some explanation of appearance of the logarithmic function at the Hardy potential in the limiting case $p=N$.

Hereinafter we set $f(t)= e^{-t}$.

%%%%%%%%%%%%%%%%%%%%%%%%%%%%%%%%%%%%%%

\noindent
{\it Proof of Theomre \ref{IL H}}：\,\,
From Lemma 2 \cite{SS}, it is enough to show the inequality (\ref{NS CH}) for radial function $u \in C_c^1(B_1)$. 
By the Hardy inequality (\ref{H}), for $k \ge 1$ we have 
\begin{align*}
\( \frac{N-p_k}{p_k} \)^{p_k} \int_{A_k} \frac{|u_k|^{p_k}}{|x|^{p_k}} dx \le \int_{A_k} | \nabla u_k |^{p_k} dx \le |A_k |^{1-\frac{p_k}{N}} \| \nabla u_k \|_N^{N-\frac{1}{k}}.
\end{align*}
From (\ref{H_p_k}) and (\ref{f}), for $k \ge 1$ we have
\begin{align*}
C \int_{A_k} \frac{|u_k|^N}{|x|^N \( \log \frac{1}{|x|} \)^N} \,dx \le \int_{A_k} | \nabla u_k |^N dx.
\end{align*}
Therefore for any $a>1$ and $k \in \Z$ we have
\begin{align}\label{H_N_k}
C \int_{A_k} \frac{|u_k|^N}{|x|^N \( \log \frac{a}{|x|} \)^\beta} \,dx \le b_k \int_{A_k} | \nabla u_k |^N dx,
\end{align}
where
\begin{align*}
b_k=
\begin{cases}
k^{N-\beta} \quad &\text{if} \,\,\,k\ge 1,\\
1 &\text{if} \,\,\,k=0, -1,\\
0 &\text{if} \,\,\,k\le -2.
\end{cases}
\end{align*}
Here, note that the inequality (\ref{H_N_k}) for $k=0, -1$ comes from the Poincar\'e inequality and the boundedness of the function $|x|^{-N}(\log \frac{a}{|x|})^{-\beta}$ on the domain $A_0 \cup A_{-1} \subset B_1 \setminus B_{e^{-2}}$. 
%Taking a sum for (\ref{H_N_k}
Since 
\begin{align*}
C \sum_{k \in \Z} \int_{B_1} \frac{|u \phi_k|^N}{|x|^N \( \log \frac{a}{|x|} \)^\beta} \,dx \le \sum_{k \in \Z} b_k \int_{A_k} | \nabla (u \phi_k ) |^N dx,
\end{align*}
we have  
\begin{align*}
C \int_{B_1} \frac{|u |^N}{|x|^N \( \log \frac{a}{|x|} \)^\beta} \,dx &\le 2^{N-1} \sum_{k \in \Z} b_k \int_{A_k} \phi_k^N | \nabla u |^N + |u|^N |\nabla \phi_k |^N dx \notag \\
&\le 2^{N-1} \int_{B_1} | \nabla u |^N dx + C \sum_{k=1}^{+\infty} b_k e^{kN} \int_{A_k} |u|^N \, dx.
\end{align*}
By Radial Lemma for $u$, we have
\begin{align*}
b_k e^{kN} \int_{A_k} |u|^N \, dx 
&\le C b_k e^{kN} \| \nabla u \|_N^N \int_{A_k} \( \log \frac{1}{|x|} \)^{N-1} \, dx \\
&\le C b_k e^{kN} \| \nabla u \|_N^N \int_k^{k+2} s^{N-1} e^{-sN} \, ds \le C b_k k^{N-1} \| \nabla u \|_N^N
\end{align*}
which implies that for $\beta > 2N$ 
\begin{align*}
C \int_{B_1} \frac{|u |^N}{|x|^N \( \log \frac{a}{|x|} \)^\beta} \,dx 
&\le C \int_{B_1} | \nabla u |^N dx + C \( \sum_{k=1}^{+\infty} k^{-1-(\beta -2N)} \) \int_{B_1} |\nabla u|^N \, dx \\
&\le C \int_{B_1} |\nabla u|^N \, dx.
\end{align*}

\qed

%%%%%%%%%%%%%%%%%%%%%%%%%%%%%%%%%%%%%%%%

\subsection{$N \nearrow \infty$}
\tc{white}{a}\\
For fixed $p\, (< N)$, we consider some infinite dimensional form of the Sobolev inequality (\ref{S}). 
Of course, since we can not consider a limit of (\ref{S}) as $N \nearrow \infty$ in the usual sense, we have to think out something, see also \cite{G} p.1062. 
In this section, we assume $p=2$. By using the scalar product strucure of the Euclidean space, we derive a logarithmic Sobolev inequality with the best constant from the Sobolev inequality (\ref{S}) as $N \nearrow \infty$ based on Beckner-Peason's paper \cite{BP}. 
We can not apply this method in $L^p$ case. For the best constant and the attainability of the $L^p$ logarithmic Sobolev inequality, we refer \cite{DD}.

\begin{theorem}\label{IL S N} (\cite{BP})
We obtain the following logarithmic Sobolev inequality (\ref{log-Sob}) as a limit of the Sobolev inequality (\ref{S}) as $N \nearrow \infty$: for any $u \in W^{1,2}(\re^n)\,\text{with}\, \int_{\re^n} |u|^2 dx =1$,
\begin{align}\label{log-Sob}
	\int_{\re^n} |u|^2 (\log |u|^2) dy \le \frac{n}{2} \log \( \frac{2}{\pi e n} \int_{\re^n} |\nabla u|^2 dy \).
\end{align}
\end{theorem}

\begin{remark}
The constant $\frac{2}{\pi e n}$ in (\ref{log-Sob}) is optimal.
\end{remark}

%%%%%%%%%%%%%%%%%%%%%%%%%%%%%%%%%%%%%%

\noindent
{\it Proof of Theorem \ref{IL S N}}：\,\,
Taking log on the both sides of the Sobolev inequality (\ref{S}) and applying Jensen's inequality, we have the logarithmic Sobolev inequality (\ref{S Jensen}) {\bf without optimal constant} for any $f \in W^{1,2}(\re^N)$ with $\int_{\re^N} |f|^2 \,dx=1$ as follows.
\begin{align}
\label{S Jensen}
\frac{2}{N} \int_{\re^N} |f|^2 (\log |f|^2) dx \le 
\frac{2}{2^*} \log \( \int_{\re^N} |f|^{2^* -2} |f|^2 \,dx \) \le \log \(  S_{N, 2}^{-1} \int_{\re^N} |\nabla f|^2 \,dx\)
\end{align}
Let $N = \ell n$ for $\ell \in \N$. By the scalar product structure of the Euclidean space, we see
\begin{align*}
	x = (\underbrace{x^1, \cdots, x^\ell}_{\ell}) \in \re^n \times \cdots \times \re^n = \re^{\ell n} = \re^N, \,\,x^i = (x^i_1, \cdots,x^i_n) \in \re^n,
\end{align*}
where $i = 1,2,\cdots,\ell$. 
For any $u \in W^{1,2}(\re^n)$ with $\int_{\re^n} |u|^2 dx =1$, set $f(x) = \prod_{i=1}^\ell u(x^i)$. Then we have 
\begin{align*}
	&\int_{\re^N} |f(x)|^2 dx = \prod_{i=1}^\ell \int_{\re^n} |u(x^i)|^2 dx^i = 1, \\
	&\int_{\re^N} |f(x)|^2 (\log |f(x)|^2 ) dx = \ell \int_{\re^n} |u(y)|^2 (\log |u(y)|^2) dy, \\
	&\int_{\re^N} |\nabla f(x)|^2 dx = \ell \int_{\re^n} |\nabla u(y)|^2 dy.
\end{align*}
Applying these equalities to (\ref{S Jensen}), we have 
\begin{align*}
\frac{2}{n} \int_{\re^n} |u|^2 (\log |u|^2) dy 
\le \log \(  \frac{1}{n \pi (N-2)} \( \frac{\Gamma (N)}{\Gamma (\frac{N}{2})} \)^{\frac{2}{N}} \int_{\re^n} |\nabla u|^2 \,dy\).
\end{align*}
Since $\Gamma (t) \sim \sqrt{2\pi} \,t^{\,t-\frac{1}{2}} \,e^{-t} \,\, \text{as}\,\, t \to \infty \,\,(\text{Stirling's formula})$, we have 
\begin{align*}
\frac{1}{n \pi (N-2)} \( \frac{\Gamma (N)}{\Gamma (\frac{N}{2})} \)^{\frac{2}{N}}
\sim \frac{2}{\pi e n}\,\, (N \to \infty)
\end{align*}
which implies the logarithmic Sobolev inequality (\ref{log-Sob}) {\bf with optimal constant}.
\qed

\begin{remark}
It is known that the Gaussian logarithmic Sobolev inequality \cite{G} is equivalent to the Euclidian logarithmic Sobolev inequality \cite{W} (see e.g. \cite{W}, \cite{B}). 
\end{remark}
%
%
%%%%%%%%%%%%%%%%%%%%%%%%%%%%%%%%%%%%%%%%%%%%%%%%%%%%%%%%%%%%%%%%
%
%\section{Direct limiting procedure}
%
%%%%%%%%%%%%%%%%%%%%%%%%%%%%%%%%%%%%%%%%%%%%%%%%%%%%%%%%%%%%%%%%
%
\section{Direct limiting procedure}

We can not take a limit directly for the Hardy inequality (\ref{H}) and the Sobolev inequality (\ref{S}) in the usual sense as $p \nearrow$ or $N \nearrow \infty$. 
In this section, we derive equivalent forms to the Hardy and Sobolev inequalities via a transformation, and we take a limit directly for these equivalent forms in the usual sense as $p \nearrow$ or $N \nearrow \infty$. For an unified viewpoint to such kind of transformations, see \cite{S(ArXiv)} Section 2. {\bf In this section, we consider only radial functions.} 

\subsection{$p \nearrow N$}
\tc{white}{a}\\
In this subsection, we refer \cite{I}. By \cite{I}, the following transformation is introduced for the Hardy and Sobolev inequalities on the whole space.
\begin{align}\label{trans}
u(r) = w(t), \,\,\text{where}\,\, r^{-\frac{N-p}{p-1}} - R^{-\frac{N-p}{p-1}} = t^{-\frac{N-p}{p-1}}
\end{align}
Here $u \in C_{{\rm rad}}^1 (B_R^N \setminus \{ 0 \}) \cap C (B_R^N), w \in C_{{\rm rad}}^1 (\re^N \setminus \{ 0 \}) \cap C(\re^N), x \in B_R^N, y \in \re^N$, $r=|x|, s=|y|$. 
Note that in the transformation (\ref{trans}), the left-hand side is the fundamental solution of $p-$Laplacian on $B_R^N$, and the right-hand side is it on $\re^N$. Since 
\begin{align*}
\frac{dr}{dt} = \(\, \frac{r}{t} \,\)^{\frac{N-1}{p-1}},
\end{align*}
we have 
\begin{align}
\int_{\re^N} | \nabla w |^p \,dy 
&= \w_{N-1} \int_0^\infty \left| \frac{d w}{d t} \,\right|^p t^{N-1} \,dt \notag \\
\label{nab}
&= \w_{N-1} \int_0^R \left| \frac{d u}{d r} \,\right|^p \( \frac{dr}{dt} \)^{p-1} t^{N-1} \,dr = \int_{B_R} | \nabla u |^p \,dx.
\end{align}
On the other hand, we have 
\begin{align}
\label{Hardy term}
\int_{\re^N} \frac{| w |^{p}}{|y|^p} \,dy &= \int_{B^N_R} \frac{| u |^{p}}{|x|^p\( 1- \(\frac{|x|}{R}\)^{\frac{N-p}{p-1}}\)^{p}} \,dx,\\
\label{Sobolev term}
\int_{\re^N} | w |^{p^*} \,dy &= \int_{B^N_R} \frac{| u |^{p^*}}{\( 1- \(\frac{|x|}{R}\)^{\frac{N-p}{p-1}}\)^{\frac{(N-1)p}{N-p}}} \,dx.
\end{align}
Therefore we see that the Sobolev inequality on the whole space for $w$:
\begin{align*}
S_{N,p} \( \int_{\re^N} |w|^{p^*} \,dy \)^{\frac{p}{p^*}} &\le \int_{\re^N} | \nabla w |^p \, dy
\end{align*}
is equivalent to the following inequality for $u$ on $B_R^N$:
\begin{align}\label{IS}
S_{N,p} \( \int_{B^N_R} \frac{| u |^{p^*}}{\( 1- \(\frac{|x|}{R}\)^{\frac{N-p}{p-1}}\)^{\frac{(N-1)p}{N-p}}} \,dx \)^{\frac{p}{p^*}} &\le \int_{B_R^N} | \nabla u |^p \, dx.
\end{align}
Since the inequality (\ref{IS}) also has the boundary singularity, we observe that the inequality (\ref{IS}) is an improvement of the classical Sobolev inequality (\ref{S}). Moreover, since the improved inequality (\ref{IS}) is equivalent to the classical Sobolev inequality (\ref{S}) on the whole space under the transformation (\ref{trans}), we can obtain several results (e.g. the scale invariance structure, the attainability of the best constant etc.) for the improved  inequality (\ref{IS}) from results for the classical Sobolev inequality (\ref{S}), see \cite{I}. 

We can not take a limit directly for the classical Sobolev inequality (\ref{S}) in the usual sense as $p \nearrow N$. 
However {\bf it is possible to take a limit directly for the improved Sobolev inequality (\ref{IS})}. Indeed, since $\lim_{x \to 0} \frac{1-r^x}{x} = \log \frac{1}{r}$ for $r \in (0,1)$, we have 
\begin{align}\label{limit}
\( 1- \(\frac{|x|}{R}\)^{\frac{N-p}{p-1}}\)^{\frac{(N-1)p}{N-p}} 
\sim \(  \frac{N-p}{p-1} \log \frac{R}{|x|} \)^{\frac{N-1}{N} p^*} \,\,(p \nearrow N).
\end{align}
Therefore, on the left-hand of the improved Sobolev inequality (\ref{IS}), we have
\begin{align*}
S_{N,p} \( \int_{B^N_R} \frac{| u |^{p^*}}{\( 1- \(\frac{|x|}{R}\)^{\frac{N-p}{p-1}}\)^{\frac{(N-1)p}{N-p}}} \,dx \)^{\frac{p}{p^*}} 
\to 
\frac{\pi^{\frac{N}{2}} N}{\Gamma(1+\frac{N}{2})} \( \sup_{x \in B_R^N} \frac{|u(x)|}{\( \log \frac{R}{|x|} \)^{\frac{N-1}{N}}} \)^N  \,\,(p \nearrow N).
\end{align*}
Hence we obtain the following.
\begin{theorem}\label{DL S}(\cite{I})
We obtain the following Alvino inequality \cite{Al} as a limit of the Sobolev inequality (\ref{S}) on the whole space as $p \nearrow N$.
\begin{align*}
\frac{\pi^{\frac{N}{2}} N}{\Gamma(1+\frac{N}{2})} \( \sup_{x \in B_R^N} \frac{|u(x)|}{\( \log \frac{R}{|x|} \)^{\frac{N-1}{N}}} \)^N \le \int_{B_R^N} | \nabla u |^N \, dx.
\end{align*}
\end{theorem}

On the other hand, we see that the Hardy inequality on the whole space for $w$ is equivalent to the following inequality on $B_R^N$ for $u$.
\begin{align}\label{IH}
\( \frac{N-p}{p} \)^p \int_{B^N_R} \frac{| u |^{p}}{|x|^p\( 1- \(\frac{|x|}{R}\)^{\frac{N-p}{p-1}}\)^{p}} \,dx
\le \int_{B_R^N} | \nabla u |^p \, dx.
\end{align}
By (\ref{limit}), we have
\begin{align*}
\( \frac{N-p}{p} \)^p \int_{B^N_R} \frac{| u |^{p}}{|x|^p\( 1- \(\frac{|x|}{R}\)^{\frac{N-p}{p-1}}\)^{p}} \,dx
\to \( \frac{N-1}{N} \)^N \int_{B^N_R} \frac{| u |^N}{|x|^N \( \log \frac{R}{|x|} \)^N} \,dx\,\,(p \nearrow N).
\end{align*}

\begin{theorem}\label{DL H}(\cite{I})
We obtain the critical Hardy inequality (\ref{CH}) as a limit of the Hardy inequality (\ref{H}) on the whole space as $p \nearrow N$.
\begin{align}\label{CH}
\( \frac{N-1}{N} \)^N \int_{B^N_R} \frac{| u |^N}{|x|^N \( \log \frac{R}{|x|} \)^N} \,dx \le \int_{B_R^N} | \nabla u |^N \, dx.
\end{align}
\end{theorem}

\begin{remark}\label{non-radial}
The above calculation holds only for radial functions. If we consider the transformation (\ref{trans}) for non-radial functions as follows:
\begin{align*}
u(r \w) = w(t \w), \,\,\text{where}\,\, r^{-\frac{N-p}{p-1} } - R^{-\frac{N-p}{p-1} } = t^{-\frac{N-p}{p-1} } \,\,\text{and}\,\,\w \in \mathbb{S}^{N-1},
\end{align*}
we have
\begin{align*}
\int_{\re^N} | \nabla w |^p \,dy = \int_{B^N_R} | L_p u |^p \,dx,
\,\,
\text{where}\,\,
L_p u = \frac{\pd u}{\pd r} \w + \frac{1}{r} \nabla_{\mathbb{S}^{N-1}} u \left[ 1 -  \( \frac{r}{R} \)^{\frac{N-p}{p-1} } \right]^{-1}.
\end{align*}
We observe that the differential operator $L_p$ is different from the usual gradient $\nabla$. 

In \cite{S(ArXiv)}, a generalization of the transformation (\ref{trans}) is considered for radial functions. Moreover, without the transformation, the attainability of minimization problems for all functions is studied in \cite{S(ArXiv)}.
\end{remark}

%%%%%%%%%%%%%%%%%%%%%%%%%%%%%%%%%%%%%%%%

\subsection{$N \nearrow \infty$}
\tc{white}{a}\\
In this subsection, we refer \cite{S(ArXiv)}. 
In order to obtain an infinite limiting form of the Sobolev inequality, we consider the following transformation:
\begin{align}\label{trans2}
u(r) = w(t), \,\,\text{where}\,\, r^{-\frac{m-p}{p-1}} = t^{-\frac{N-p}{p-1}}.
\end{align}
Here $u \in C_{{\rm rad}}^1 (\re^m \setminus \{ 0 \}) \cap C (\re^m), w \in C_{{\rm rad}}^1 (\re^N \setminus \{ 0 \}) \cap C(\re^N), x \in \re^m, y \in \re^N$, $r=|x|, s=|y|$.  
Let $p, m, N$ satisfy $1 \le p < m \le N$. 
Thanks to the transformation (\ref{trans2}), we can obtain an equivalent inequality on the lower dimensional Soboelv space $W_0^{1,p}(\re^m)$ to the Sobolev inequality on the higher dimensional Sobolev space $W_0^{1,p}(\re^N)$. Therefore, since we can regard the dimension $N$ as a parameter, we can take a limit of the Sobolev inequality as $N \nearrow \infty$. In the same way as before, we have
\begin{align*}
\frac{dr}{dt} = \frac{N-p}{m-p} \(\, \frac{r^{m-1}}{t^{N-1}} \,\)^{\frac{1}{p-1}}.
\end{align*}
Therefore we have
\begin{align*}
\int_{\re^N} | \nabla w |^p \,dy 
&= \frac{\w_{N-1}}{\w_{m-1}} \( \frac{N-p}{m-p} \)^{p-1} \int_{\re^m} | \nabla u |^p \,dx,\\
\int_{\re^N} | w |^{\frac{Np}{N-p}} \,dy 
&= \frac{\w_{N-1}}{\w_{m-1}} \,\frac{m-p}{N-p} \int_{\re^m} \frac{| u |^{\frac{Np}{N-p}}}{|x|^{\frac{N-m}{N-p}p}} \,dx.
\end{align*}
Thus we observe that the higher dimensional Sobolev inequality (\ref{S}) for $w$ on $\re^N$ is equivalent to the following lower dimensional inequality (\ref{S another}) for $u$ on $\re^m$.
\begin{align}\label{S another}
S_{N,p} \( \frac{\w_{m-1}}{\w_{N-1}} \)^{\frac{p}{N}} \( \frac{m-p}{N-p} \)^{p-\frac{p}{N}} \( \int_{\re^m} \frac{| u |^{\frac{Np}{N-p}}}{|x|^{\frac{N-m}{N-p}p}} \,dx \)^{\frac{N-p}{N}} \le \int_{\re^m} | \nabla u |^p \, dx.
\end{align}
Since
\begin{align*}
\w_{N-1} = \frac{N \pi^{\frac{N}{2}}}{\Gamma \( 1+ \frac{N}{2} \)} \,\,\text{and}\,\, 
\Gamma (t) \sim \sqrt{2\pi} \,t^{\,t-\frac{1}{2}} \,e^{-t} \,\, \text{as}\,\, t \to \infty \,\,(\text{Stirling's formula}),
\end{align*}
we can calculate the limit of the coefficient on the left-hand side of the inequality (\ref{S another}) as $N \nearrow \infty$ as follows.
\begin{align*}
&S_{N,p} \( \frac{\w_{m-1}}{\w_{N-1}} \)^{\frac{p}{N}} \( \frac{m-p}{N-p} \)^{p-\frac{p}{N}} \\
&= \pi^{\frac{p}{2}} N \( \frac{N-p}{p-1} \)^{p-1} \( \frac{m-p}{N-p} \)^{p}  \( \frac{\Gamma (\frac{N}{p}) \Gamma (N+ 1 -\frac{N}{p}) \,\w_{m-1} \,(N-p)}{\Gamma (N)  \Gamma(1+\frac{N}{2}) \, \w_{N-1} \,(m-p) } \)^{\frac{p}{N}} \\
&= \frac{N}{N-p} \,\frac{(m-p)^p}{(p-1)^{p-1}} \( \frac{\w_{m-1} \,(N-p)}{m-p} \)^{\frac{p}{N}}  \( \frac{\Gamma (\frac{N}{p}) \Gamma \( \frac{p-1}{p} N +1 \)}{\Gamma (N+1) } \)^{\frac{p}{N}} \\
&\sim \frac{(m-p)^p}{(p-1)^{p-1}} \( \frac{(\frac{N}{p})^{\frac{N}{p} -\frac{1}{2}} e^{-\frac{N}{p}} \( \frac{p-1}{p} N +1\)^{\frac{p-1}{p} N +\frac{1}{2}} e^{-\frac{p-1}{p} N -1}}{(N+1)^{N+\frac{1}{2}} e^{-(N+1)} } \)^{\frac{p}{N}} \sim \( \frac{m-p}{p} \)^p \,\,(N \nearrow \infty).
\end{align*}
Therefore, on the left-hand side of the inequality (\ref{S another}), we have
\begin{align*}
S_{N,p} \( \frac{\w_{m-1}}{\w_{N-1}} \)^{\frac{p}{N}} \( \frac{m-p}{N-p} \)^{p-\frac{p}{N}} &\( \int_{\re^m} \frac{| u |^{\frac{Np}{N-p}}}{|x|^{\frac{N-m}{N-p}p}} \,dx \)^{\frac{N-p}{N}} \to
\( \dfrac{m-p}{p} \)^p \int_{\re^m} \dfrac{|u(x)|^p}{|x|^p} dx.
\end{align*}
Hence we observe a new relationship between the Hardy and Sobolev inequalities as follows.

\begin{theorem}
We obtain the Hardy inequality (\ref{H}) as an infinite dimensional form  of the Sobolev inequality (\ref{S}).
\end{theorem}

%%%%%%%%%%%%%%%%%%%%%%%%%%%%%%%%%%%%%%%%

On the other hand, under (\ref{trans2}), we have 
\begin{align*}
\( \dfrac{N-p}{p} \)^p \int_{\re^N} \frac{| w |^p}{|y|^p} \,dy 
= \frac{\w_{N-1}}{\w_{m-1}} \( \frac{N-p}{m-p} \)^{p-1} \( \dfrac{m-p}{p} \)^p \int_{\re^m} \frac{| u |^p}{|x|^p} \,dx.
\end{align*}
Therefore we see that the higher dimensional Hardy inequality for $w$ on $\re^N$ is equivalent to the lower dimensional Hardy inequality for $u$ on $\re^m$. 
Hence, we observe that {\bf the Hardy inequality (\ref{H}) is dimension free in some sence}, and an infinite dimensional form of the Hardy inequality is the Hardy inequality again.

%%%%%%%%%%%%%%%%%%%%%%%%%%%%%%%%%%%%%%%%%%%%%%%%%%%%%%%%%%%%%%%%%%%%%%%%%%%%%%%%%%%
%
% \section{Summary and supplement}
%
%%%%%%%%%%%%%%%%%%%%%%%%%%%%%%%%%%%%%%%%%%%%%%%%%%%%%%%%%%%%%%%%%%%%%%%%%%%%%%%%%%%
\section{Summary and supplement}
\noindent
We summarize \S 2 and \S 3.

In \S 2, we consider two limits of the Hardy and Sobolev inequalities as $p \nearrow N$ and $N \nearrow \infty$ via indirect limiting procedures. For these indirect limiting procedures, we can not expect to get an information of the best constant in the limiting inequality in general. 
However it is possible to apply the indirect limiting procedure to a higher order inequality and another inequality. 
Indeed, in \cite{SS}, the indirect limiting procedure is applied to the Rellich inequality, which is known as a higher order generalization of the Hardy inequality, and the Poincar\'e inequality (For the Rellich inequality, we consider a limit as $p \nearrow \frac{N}{2}$. For the Poincar\'e inequality, we consider a limit as $|\Omega| \searrow 0$).

On the other hand, in \S 3, we restrict radial functions only and we derive equivalent forms to the Hardy and Sobolev inequalities via some transformations. Through these equivalent forms, we consider two limits of the Hardy and Sobolev inequality as $p \nearrow N$ and $N \nearrow \infty$ via direct limiting procedures. For these direct limiting procedures, we can obtain an information of the best constant in the limiting inequality. However, these direct limiting procedures are based on the special transformations. Therefore it seems difficult to generalize to a higher order case. 
Based on these transformations, we observe that two well-known embeddings of the subcritical Sobolev space $W_0^{1,p} (p<N)$：
\begin{align*}
&W_0^{1, p} \hookrightarrow L^{p^*, p} \,(\text{The Hardy inequality}), \\
&W_0^{1, p} \hookrightarrow L^{p^*, p^*}= L^{p^*} \,(\text{The Sobolev inequality})
\end{align*} 
become the following embeddings in the limiting case where $p=N$.
\begin{align*}
&W_0^{1, N} (B_1) \hookrightarrow L^{\infty, N}(\log L)^{-1} \,(\text{The critical Hardy inequality}),\\
&W_0^{1, N} (B_1) \hookrightarrow L^{\infty, \infty}(\log L)^{-1+\frac{1}{N}} = {\rm Exp L}^{\frac{N}{N-1}}\\
&\hspace{5em}(\text{The Alvino inequality, The Trudinger-Moser inequality}).
\end{align*} 
From an inclusion property of the Lorentz-Zygmund space (see e.g. \cite{BS} Theorem 9.5.), we obtain the embeddings of the critical Sobolev space $W_0^{1,N} (B_1)$ as follows：
\begin{align*}
W_0^{1,N} (B_1) \hookrightarrow L^{\infty, N}(\log L)^{-1} \hookrightarrow L^{\infty, q}(\log L)^{-1+ \frac{1}{N} - \frac{1}{q}} \hookrightarrow L^{\infty, \infty}(\log L)^{-1+\frac{1}{N}} = {\rm Exp L}^{\frac{N}{N-1}}
\end{align*}
for any $q \in (N, \infty)$. For the attainability of the best constants associated with the above embeddings (inequalities), see \cite{Al, CC, AS, HK, CRT(2013), II, S(JDE)}.

In addition, we refer several related works to this note as follows: 
\begin{itemize}
\item \cite{Y}： $L^p$ boundedness of the Hilbert transformation when  $p \searrow 1$
\item \cite{Z} XII 4.41.： $L^p$ boundedness of the Hilbert transformation when $p \nearrow \infty$
\item \cite{SC} Corollary 3.2.4： A derivation of the Sobolev inequality form the Nash inequality
\item \cite{F}： The $L^p$ logarithmic Sobolev inequality when $p \nearrow \infty$
\item \cite{CRT(2018)}： A derivation of some {\it equivalent} inequality to the classical Hardy inequality (Here, the meaning of the equivalence in \cite{CRT(2018)} is weaker than it in \S 3. That is, the validity of two equivalent inequalities is corresponding each other, but the attainability of two best constants in two equivalent inequalities is not corresponding each other.)
\end{itemize}

\vspace{1em}

%%%%%%%%%%%%%%%%%%%%%%%%%%%%%%%%%%%%%%%%%%%%%%%%%%%%%%%%%%%%%%%%%%%%%%%%%%%%%%%%%%%
%
% \section{Acknowledgement}
%
%%%%%%%%%%%%%%%%%%%%%%%%%%%%%%%%%%%%%%%%%%%%%%%%%%%%%%%%%%%%%%%%%%%%%%%%%%%%%%%%%%%
\noindent
{\large {\bf Acknowledgement}}

This work was (partly) supported by Osaka City University Advanced
Mathematical Institute (MEXT Joint Usage/Research Center on Mathematics
and Theoretical Physics). 
The author was supported by JSPS KAKENHI Early-Career Scientists, No. JP19K14568.

\end{document}